\documentclass{amsart}

\newtheorem{theorem}{Theorem}[section]
\newtheorem{lemma}[theorem]{Lemma}

\theoremstyle{definition}

\theoremstyle{remark}

\numberwithin{equation}{section}

\begin{document}

\title{Energy and volume of vector fields on spherical domains}

\author{Fabiano G. B. Brito}
\address{Universidade Federal do Estado do Rio de Janeiro, Brazil}

\author{Andr\'e Gomes}
\address{Universidade de S\~ao Paulo, Brazil}
\author{Giovanni S. Nunes}
\address{Universidade Federal de Pelotas, Brazil}




\begin{abstract}
We present in this paper a ``boundary version" for theorems about minimality of volume and energy functionals on a spherical domain of three-dimensional Euclidean sphere. 
\end{abstract}

\maketitle

\section{Introduction}
Let $(M,g)$ be a closed, n-dimensional Riemannian manifold and $T^{1}M$ the unit tangent bundle of $M$ considered as a closed Riemannian manifold with the Sasaki
metric. Let $X:M\longrightarrow T^{1}M$ be a unit vector field defined on $M$, regarded as a smooth section on the unit tangent bundle $T^{1}M$. The volume of $X$ was defined in [8] by $\mathrm{vol}(X):=\mathrm{vol}(X(M))$, where $\mathrm{vol}(X(M))$ is the volume of the submanifold $X(M)\subset T^{1}M$. Using an orthonormal local frame $\left\{e_{1}, e_{2},\ldots,e_{n-1}, e_{n}=X\right\}$, the volume of the unit vector field $X$ is given by
\begin{eqnarray*}
\mathrm{vol}(X)=\int_{M} (1+\sum\limits_{a=1}^{n}\left\|\nabla_{e_{a}}X\right\|^{2}+\sum\limits_{a<b}\left\|\nabla_{e_{a}}X\wedge\nabla_{e_{b}}X\right\|^{2}+\ldots
\\
\ldots+\sum\limits_{a_{1}<\cdots<a_{n-1}}\left\|\nabla_{e_{a_{1}}}X\wedge\cdots\wedge\nabla_{e_{a_{n-1}}}X\right\|^{2})^{1/2} \nu_{_{M}}(g)	
\end{eqnarray*}
and the energy of the vector field $X$ is given by 
\begin{eqnarray*}
	\mathcal{E}(X)=\frac{n}{2}\mathrm{vol}(M)+\frac{1}{2}\int_{M}\sum\limits_{a=1}^{n}\left\|\nabla_{e_{a}}X\right\|^{2}\nu_{_{M}}(g)
\end{eqnarray*}
The Hopf vector fields on $\mathbb{S}^{3}$ are unit vector fields tangent to the classical Hopf
fibration $\pi:\mathbb{S}^{3}\longrightarrow \mathbb{S}^{2}$ with fiber homeomorphic to $\mathbb{S}^{1}$. 
\\
The following theorems gives a characterization of Hopf flows as absolute minima of volume and energy functionals:

\theorem [\textbf{[8]}] {The unit vector fields of minimum volume on the sphere $\mathbb{S}^{3}$ are precisely
the Hopf vector fields and no others.}

\normalfont
\theorem [\textbf{[1]}] {The unit vector fields of minimum energy on the sphere $\mathbb{S}^{3}$ are precisely
the Hopf vector fields and no others.}
\\
\\
\normalfont We prove in this paper the following boundary version for these Theorems:

\theorem {Let $U$ be an open set of the three-dimensional unit sphere $\mathbb{S}^{3}$ and let $K\subset U$ be a compact set. Let $\vec{v}$ be an unit vector field on $U$ which coincides with a Hopf flow $H$ along the boundary of K. Then $\mathrm{vol}(\vec{v})\geq \mathrm{vol}(H)$ and $\mathcal{E}(\vec{v})\geq \mathcal{E}(H)$.}
\\
\\
\normalfont
Other results for higher dimensions may be found in [2], [5], [7] and [8].

\section{Preliminaries}
Let $U\subset\mathbb{S}^{3}$ be an open set. We consider a compact set $K\subset U$. Let $H$ be a Hopf vector field on $\mathbb{S}^{3}$ and let $\vec{v}$ be an unit vector field defined on $U$. We also consider the map $\varphi_{t}^{\vec{v}}:U\longrightarrow \mathbb{S}^{3}(\sqrt{1+t^{2}})$ given by $\varphi_{t}^{\vec{v}}(x)=x+t\vec{v}(x)$. This map was introduced in [10] and [3].

\lemma {For $t>0$ sufficiently small, the map $\varphi_{t}^{\vec{v}}$ is a diffeomorphism.}
\normalfont

\proof A simple application of the identity perturbation method $\square$
\\
\\
In order to find the Jacobian matrix of $\varphi_{t}^{\vec{u}}$, we define the unit vector field $\vec{u}$ 
\begin{eqnarray*}
	\vec{u}(x):=\frac{1}{\sqrt{1+t^{2}}}\vec{v}(x)-\frac{t}{\sqrt{1+t^{2}}}x
\end{eqnarray*}
Using an adapted orthonormal frame $\left\{e_{1},e_{2},\vec{v}\right\}$ on a neighborhood $V\subset U$, we obtain an adapted orthonormal frame on $\varphi_{t}^{\vec{v}}(V)$ given by $\left\{\bar{e}_{1},\bar{e}_{2},\vec{u}\right\}$, where $\bar{e}_{1}=e_{1}$, $\bar{e}_{2}=e_{2}$.
\\
\\
In this manner, we can write
\begin{eqnarray*}
	d\varphi_{t}^{\vec{v}}(e_{1})\!\!\!&=&\!\!\!\left\langle d\varphi_{t}^{\vec{v}}(e_{1}),e_{1}\right\rangle e_{1}+\left\langle d\varphi_{t}^{\vec{v}}(e_{1}),e_{2}\right\rangle e_{2}+\left\langle d\varphi_{t}^{\vec{v}}(e_{1}),\vec{u}\right\rangle \vec{u}\\
	d\varphi_{t}^{\vec{v}}(e_{2})\!\!\!&=&\!\!\!\left\langle d\varphi_{t}^{\vec{v}}(e_{2}),e_{1}\right\rangle e_{1}+\left\langle d\varphi_{t}^{\vec{v}}(e_{2}),e_{2}\right\rangle e_{2}+\left\langle d\varphi_{t}^{\vec{v}}(e_{2}),\vec{u}\right\rangle \vec{u}\\
	d\varphi_{t}^{\vec{v}}(\vec{v})\!\!\!&=&\!\!\!\left\langle d\varphi_{t}^{\vec{v}}(\vec{v}),e_{1}\right\rangle e_{1}+\left\langle d\varphi_{t}^{\vec{v}}(\vec{v}),e_{2}\right\rangle e_{2}+\left\langle d\varphi_{t}^{\vec{v}}(\vec{v}),\vec{u}\right\rangle \vec{u}
\end{eqnarray*}
Now, by Gauss' equation of immersion $\mathbb{S}^{3}\hookrightarrow \mathbb{R}^{4}$, we have 
\begin{eqnarray*}
	d\vec{v}(Y)=\nabla_{Y}\vec{v}-\left\langle \vec{v},Y\right\rangle x
\end{eqnarray*}
for every vector field $Y$ on $\mathbb{S}^{3}$, and then
\begin{eqnarray*}
	\left\langle d\varphi_{t}^{\vec{v}}(e_{1}),e_{1}\right\rangle = \left\langle e_{1}+td\vec{v}(e_{1}),e_{1}\right\rangle = 1 + t\left\langle \nabla_{e_{1}}\vec{v},e_{1}\right\rangle
\end{eqnarray*}
Analogously, we can conclude that 
\begin{eqnarray*}
\left\langle d\varphi_{t}^{\vec{v}}(e_{1}),e_{2}\right\rangle \!\!\!&=&\!\!\!t\left\langle \nabla_{e_{1}}\vec{v},e_{2}\right\rangle\\
\left\langle d\varphi_{t}^{\vec{v}}(e_{2}),e_{1}\right\rangle \!\!\!&=&\!\!\! t\left\langle \nabla_{e_{2}}\vec{v},e_{1}\right\rangle\\
\left\langle d\varphi_{t}^{\vec{v}}(e_{2}),e_{2}\right\rangle \!\!\!&=&\!\!\! 1+t\left\langle \nabla_{e_{2}}\vec{v}, e_{2}\right\rangle\\
\left\langle d\varphi_{t}^{\vec{v}}(e_{1}),\vec{u}\right\rangle \!\!\!&=&\!\!\! 0\\
\left\langle d\varphi_{t}^{\vec{v}}(e_{2}),\vec{u}\right\rangle \!\!\!&=&\!\!\! 0\\
\left\langle d\varphi_{t}^{\vec{v}}(\vec{v}),\vec{u}\right\rangle \!\!\!&=&\!\!\! \sqrt{1+t^{2}}
\end{eqnarray*}
By applying the notation $h_{ij}(\vec{v}):=\left\langle \nabla_{e_{i}}\vec{v},e_{j}\right\rangle$ ($i,j=1,2$), the determinant of the Jacobian matrix of $\varphi_{t}^{\vec{v}}$ can be express in the form 
\begin{eqnarray*}
	\det(d\varphi_{t}^{\vec{v}})=\sqrt{1+t^{2}}(1+\sigma_{1}(\vec{v}).t+\sigma_{2}(\vec{v}).t^{2})
\end{eqnarray*}
where, by definition, 
\begin{eqnarray*}
\sigma_{1}(\vec{v})\!\!\!&:=&\!\!\!h_{11}(\vec{v})+h_{22}(\vec{v})\\ \sigma_{2}(\vec{v})\!\!\!&:=&\!\!\!h_{11}(\vec{v})h_{22}(\vec{v})-h_{12}(\vec{v})h_{21}(\vec{v})
\end{eqnarray*}

\section{Proof of Theorem 1.3}

The energy of the vector field $\vec{v}$ (on $K$) is given by 
\begin{eqnarray*}	\mathcal{E}(\vec{v}):=\frac{1}{2}\int_{K}\left\|d\vec{v}\right\|^{2}=\frac{3}{2}\mathrm{vol}(K)+\frac{1}{2}\int_{K}\left\|\nabla \vec{v}\right\|^{2}
\end{eqnarray*}

Using the notations above, we have
\begin{eqnarray*}	\mathcal{E}(\vec{v})\!\!\!&=&\!\!\!\frac{3}{2}\mathrm{vol}(K)+\frac{1}{2}\int_{K}[(\sum\limits_{i,j=1}^{2}(h_{ij})^{2}+(\left\langle \nabla_{\vec{v}}\vec{v},e_{1}\right\rangle)^{2}+(\left\langle \nabla_{\vec{v}}\vec{v},e_{2}\right\rangle)^{2}]
\end{eqnarray*}
and then
\begin{eqnarray*}
\mathcal{E}(\vec{v})&\geq&\!\!\! \frac{3}{2}\mathrm{vol}(K)+\frac{1}{2}\int_{K}\sum\limits_{i,j=1}^{2}(h_{ij})^{2}\\
\!\!\!&\geq&\!\!\! \frac{3}{2}\mathrm{vol}(K)+\frac{1}{2}\int_{K}2(h_{11}h_{22}-h_{12}h_{21})\\
\!\!\!&=&\!\!\!\frac{3}{2}\mathrm{vol}(K)+\int_{K}\sigma_{2}(\vec{v})
\end{eqnarray*}

On the other hand, by change of variables theorem, we obtain
\begin{eqnarray*}	\mathrm{vol}[\varphi_{t}^{H}(K)]=\int_{K}\sqrt{1+t^{2}}(1+\sigma_{1}(H).t+\sigma_{2}(H).t^{2})=\delta\cdot\mathrm{vol}(\mathbb{S}^{3}(\sqrt{1+t^{2}}))
\end{eqnarray*}
where $\delta:=\mathrm{vol}(K)/\mathrm{vol}(\mathbb{S}^{3})$.
\\
(Remark that $\sigma_{1}(H)$ and $\sigma_{2}(H)$ are constant functions on $\mathbb{S}^{3}$, in fact, we have $\sigma_{1}(H)=0$ and $\sigma_{2}(H)=1$, by a straightforward computation shown in [6]).
\\
\\
Suppose now that $\vec{v}$ is an unit vector field on $K$ which coincides with a Hopf vector field $H$ on the boundary of $K$. Then, obviously
\begin{eqnarray*}
	\mathrm{vol}[\varphi_{t}^{\vec{v}}(K)]=\mathrm{vol}[\varphi_{t}^{H}(K)]
\end{eqnarray*}

Therefore, we obtain 
\begin{eqnarray*}	\mathrm{vol}[\varphi_{t}^{\vec{v}}(K)]\!\!\!&=&\!\!\!\int_{K}\sqrt{1+t^{2}}(1+\sigma_{1}(\vec{v}).t+\sigma_{2}(\vec{v}).t^{2})\\
&=&\!\!\!\delta\cdot\mathrm{vol}(\mathbb{S}^{3}(\sqrt{1+t^{2}}))=[\mathrm{vol}(K)](1+t^{2})^{3/2}
\end{eqnarray*}
By identity of polynomials, we conclude that 
\begin{eqnarray*}
	\int_{K}\sigma_{2}(\vec{v})=\mathrm{vol}(K)
\end{eqnarray*}
and consequently
\begin{eqnarray*}
	\mathcal{E}(\vec{v})\geq \frac{3}{2}\mathrm{vol}(K)+\mathrm{vol}(K)=\mathcal{E}(H)
\end{eqnarray*}
Now, observing that 
\begin{eqnarray*}
\mathrm{vol}(H)= 2\mathrm{vol}(K), \ \ \ \ \int_{K}\sigma_{2}(\vec{v})=\mathrm{vol}(K)\ \ \ and \ \ \ \sum\limits_{i,j=1}^{2} h_{ij}^{2}(\vec{v})\geq 2\sigma_{2}(\vec{v})
\end{eqnarray*}
we can obtain an analogue of this result for volumes
\begin{eqnarray*}
\mathrm{vol}(X)\!\!\!&=&\!\!\!\int_{K}\sqrt{1+\sum\limits h_{ij}^{2}+[\det(h_{ij})]^{2}+\cdots}
\\
&\geq&\!\!\!\int_{K}\sqrt{1+2\sigma_{2}+\sigma_{2}^{2}}
\\
&=&\!\!\!\int_{K}(1+\sigma_{2})=2\mathrm{vol}(K)=\mathrm{vol}(H)\ \square
\end{eqnarray*}

\section{Final remarks}
\begin{enumerate}
  \item If $K$ is a spherical cap (the closure of a connected open set
with round boundary of the three unit sphere), the theorem provides a
``boundary version" for the minimalization theorem of
energy and volume functionals on [1] and [8].
\\
	\item The ``Hopf boundary" hypothesis is essential. In fact, if there is no constraint for the unit vector field $\vec{v}$ on $\partial K$, it is possible to construct vector fields on ``small caps" such that $\left\|\nabla \vec{v}\right\|$ is small on $K$ (exponential maps may be used on that construction). A consequence of this is that $\mathcal{E}(\vec{v})$ and $\mathrm{vol}(\vec{v})$ are less than volume and energy of Hopf vector fields respectively.
	\\
	\item The results of this paper may, possibly, be extended for the energy of solenoidal unit vector fields in the higher dimensional case ($n=2k+1$). We intend to treat this subject in a forthcoming paper.
	\\
	\item We express our gratitude to Prof. Jaime Ripoll for helpful
conversation concerning the final draft of our paper.
\end{enumerate}


\begin{thebibliography}{}

\bibitem[1]{01}F. G. B. Brito, \textit{Total bending of flows with mean curvature correction}, Diff. Geom. Appl. 12, (2000), 157-163.
\bibitem[2]{02}F. B. B. Brito, P. M. Chac\'on and A. M. Naveira, \textit{On the volume of unit vector fields on spaces of constant sectional curvature}, Comment. Math. Helv., 79 (2004) 300-316.
\bibitem[3]{03}F. G. B. Brito, R. Langevin and H. Rosenberg, \textit{Int\'egrales de courbure sur des vari\'et\'es feuillet\'ees}, J. Differential Geom. 16, (1981), no. 1, 19?50.
\bibitem[4]{04}E. Boeckx, J. C. Gonz\'ales-D\'avila and L. Vanhecke, \textit{Energy of radial vector fields on compact rank one symmetric spaces}, Ann. Global Anal. Geom., 23 (2003), no. 1, 29-52. 
\bibitem[5]{05}V. Borrelli and O. Gil-Medrano, \textit{A critical radius for unit Hopf vector fields on spheres}, Math. Ann. 334, (2006), no. 4, 731?751.
\bibitem[6]{06}P. M. Chac\'on, \textit{Sobre a energia e energia corrigida de campos unit\'arios e distribui\c{c}\~oes. Volume de campos unit\'arios}, PhD Thesis, Universidade de S\~ao Paulo, Brazil. 
\bibitem[7]{07}P. M. Chac\'on, A. M. Naveira and J. M. Weston, \textit{On the energy of distributions, with application to the quaternionic Hopf fibrations}, Monatsh. Math. 133, (2001), no. 4, 281-294. 
\bibitem[8]{08}H. Gluck and W. Ziller, \textit{On the volume of the unit vector fields on the three sphere}, Comment. Math. Helv. 61, (1986), 177-192.
\bibitem[9]{09}D. L. Johnson, \textit{Volumes of flows}, Proc. Amer. Math. Soc. 104 (1988), no. 3, 923?931.
\bibitem[10]{10}J. Milnor, \textit{Analytic proofs of the ``hairy ball theorem" and the Brouwer fixed-point theorem}, Amer. Math. Monthly 85, (1978), no. 7, 521?524.

\end{thebibliography}
\end{document}